\documentclass[11pt]{article}
\usepackage{latexsym}
\title{The Ricci flow on the two ball with a rotationally symmetric metric} 
\author{Jean C. Cortissoz}
\date{}

\newtheorem{thm}{Theorem}[section]
\newtheorem{prop}[thm]{Proposition}
\newtheorem{lem}[thm]{Lemma}

\newtheorem{theorem}{Theorem}

\newtheorem{defn}[thm]{Definition}

\newtheorem{rem}[thm]{Remark}
 
\begin{document}
\maketitle
\begin{abstract}
 In this paper we study a boundary value problem for the Ricci flow in the
 two-dimensional ball endowed with a rotationally symmetric metric of positive 
 Gaussian curvature. We show
 short and long time existence results. We construct families of metrics
 for which the flow uniformizes the curvature along a sequence of times. Finally, we show that if the
 initial metric has positive Gaussian curvature and the boundary has positive
 geodesic curvature then the flow uniformizes the curvature along a sequence of times. 

\end{abstract}


\section{Introduction.}\label{intro}
 
Let $\mathbf{B}^2$ be the two dimensional ball. We let 
$g$ be a time dependent family of metrics. Let $R$ denote its scalar curvature
and $k_g$ the geodesic curvature of its boundary $S^1 \approx \partial \mathbf{B}^2$ with respect to the
outward unit normal. 
In this paper we begin the study the of following boundary value problem for the Ricci flow
in surfaces,
\begin{equation}
\label{BVP1}
\left\{
\begin{array}{l}
\frac{\partial g}{\partial t}=-Rg \quad \mbox{in} \quad M\times\left(0,T\right), \\
k_g=\psi \quad \mbox{on} \quad \partial M \times \left(0,T\right),\\
g\left(\cdot,0\right)=g_0 \quad .
\end{array}
\right. 
\end{equation}

In other words we want to study the boundary value problem consisting of prescribing the
geodesic curvature of the boundary.

We also consider the normalized ``version'' of (\ref{BVP1}). The normalized Ricci flow is obtained from
(\ref{BVP0}) by the following procedure:

Let $\phi\left(t\right)$ be such that for $\tilde{g}=\phi g$ we have $Vol\left(M\right)=1$. Then
we change the time scale by letting
\[
\tilde{t}\left(t\right)=\int_0^t \phi\left(\tau\right)\,d\tau ,
\]
then $\tilde{g}\left(\tilde{t}\right)$ satisfies the equation
\begin{equation}
\label{normalizedflow}
\frac{\partial}{\partial \tilde{t}}\tilde{g}=\frac{2}{n}\tilde{r}\tilde{g}-2Ric\left(\tilde{g}\right) .
\end{equation}

Flows in surfaces with boundary have been studied by Brendle. In the paper \cite{Bre}, the following result is shown,
\begin{theorem}
Let $M$ be a compact surface with boundary $\partial M$. Then for every initial metric with vanishing geodesic curvature
at the boundary, the initial value problem (\ref{BVP1}) has a unique solution. The solution to the normalized flow
is defined for all $t\geq 0$. For $t\rightarrow \infty$ the solution converges exponentially to a metric
with constant Gauss curvature and vanishing geodesic curvature.
\end{theorem} 

Due to its simplicity and non triviality, and in order to give a strong basis to
our general intuition, in this paper we mainly consider rotationally symmetric metrics on $\mathbf{B}^2$
(but we should mention that the results of Section \ref{Perelman} are also valid without this restriction), i.e.
metrics of the form
\[
ds^2=dr^2+f\left(r\right)^2d\omega^2 ,
\]
and the following problem boundary value problem,
\begin{equation}
\label{BVP0}
\left\{
\begin{array}{l}
\frac{\partial g}{\partial t}=-Rg \quad \mbox{in} \quad B^2\times \left(0,T\right),\\
k_g=k_0 \quad \mbox{on} \quad \partial B^2 \times \left(0,T\right),\\
g\left(\cdot,0\right)=g_0 \quad ,
\end{array}
\right.
\end{equation} 
i.e., we fix $\psi=k_0$, where $k_0$ is the geodesic curvature of the boundary with respect to the metric $g_0$. Notice that
if $g_0$ is rotationally symmetric, then the solution $g$ of (\ref{BVP0}) will also be rotationally symmetric.

In this paper we prove the following results 
\begin{theorem}
\label{longtimeexistence}
Assume $R\geq 0$ at $t=0$. Then the unique solution to the normalized flow exists for all time.
\end{theorem}

In Section \ref{Perelman}  by using the techniques introduced by Perelman \cite{Per1}, we show the following
result.

\begin{theorem}
\label{convexboundarycase}
If $R>0$ and $k_0\geq 0$ for the initial metric, then the curvature $R$ of the solution
of (\ref{BVP0}) blows up in finite time. Also, there is a sequence of times $t_k\rightarrow T$ so that for
the solution of (\ref{BVP0})
\[
\frac{R_{max}\left(t_k\right)}{R_{min}\left(t_k\right)}\rightarrow 1
\quad \mbox{as}\quad t_k\rightarrow T .
\]
\end{theorem}

Finally, in Section \ref{FamiliesOfExamples}, we show some
families of examples where $k_0<0$, and there is a sequence of times $t\rightarrow \infty$ such that
$g\left(t\right)\rightarrow g\left(\infty\right)$ smoothly and the metric $g\left(\infty\right)$
has constant curvature and totally geodesic boundary.

The main tool we use to prove Theorem \ref{convexboundarycase}
(and any other convergence result in this paper) is blow up analysis. In order to produce blow up limits we need
to bound derivatives of the curvature in terms of a bound on the curvature. This is the content of Section \ref{DerivativeEstimates}.
Injectivity radius estimates and Fermi inradius estimates (``how big a collar of the boundary can be'')
are easy to get in our situation thanks to the symmetries assumed and by means of the Laplacian Comparison
Theorem. 
 
We must point out that the for the normalization of (\ref{BVP0}), $k_g\rightarrow 0$ (which is
due to the fact that we are rescaling by a factor that is going to $\infty$, and we keep the
geodesic curvature in the unnormalized flow uniformly bounded), and hence the normalized Ricci flow is 
uniformizing the curvature while the boundary of the surface is becoming 
totally geodesic. 

Finally,
We want to remark that with a little more work we can show that the subsequential
convergence is actually uniform convergence. Also, one can show that the
convergence of the curvature to a constant is exponential in any $L^p$ norm
for $1\leq p< \infty$.


\section{Short time existence.}\label{ShortTimeExistence}

From (\ref{BVP1}) follows that the solution $g$ is in the same conformal class
of $g_0$. If we let $\tilde{g}=ug$, for $u>0$, it is well known that 
\begin{equation}
\tilde{R}=\frac{1}{u}\left(R-\Delta \log u\right)
\end{equation} 

Let $\frac{\partial}{\partial \eta}$ be the outward unit normal to $\partial M$ with respect to the
initial
metric $g_0$. Then the transformation law for the geodesic curvature is given by
\begin{equation}
k_{\tilde{g}}=\frac{1}{\sqrt{u}}\left(k_g+\frac{1}{2u}\frac{\partial u}{\partial \eta}\right).
\end{equation}

Therefore, equation (\ref{BVP1}) is equivalent to the following Boundary Value Problem
\begin{equation}
\left\{
\label{BVP2}
\begin{array}{l}
u_t=\Delta \log u -R_0 \quad \mbox{in}\quad M\times\left(0,T\right)\\
u_{\eta}= 2u\left(k_g u^{\frac{1}{2}}-k_0\right) \quad \mbox{on} \quad \partial M \times \left(0,T\right)\\
u\left(x,0\right)=1 \quad .
\end{array}
\right.
\end{equation}

From (\ref{BVP2}), by the theory of linear parabolic equations, 
 an application of the Inverse Function Theorem, gives us the following
\begin{thm}[Short time existence]
Equation (\ref{BVP2}) has a unique solution for a short time. This solution
is smooth in $\overline{M}\times\left(0,\epsilon\right) $ and is of class $C^2$
in $\overline{M}\times\left[0,\epsilon\right)$
\end{thm}

Also, we can assert the following fact concerning regularity
\begin{thm}
If $g_0$, the initial metric, is smooth then the solution to (\ref{BVP1}) is 
smooth in $\overline{M}\times\left(0,\epsilon\right)$.
\end{thm}


\section{Evolution of the curvature.}\label{EvolutionOfTheCurvature}

It is important to know how the curvature quantities are evolving if the metric is evolving by the Ricci flow.
So, if $g$ evolves by the unnormalized Ricci flow in the time interval 
$\left(0,T\right)$, and $\frac{\partial}{\partial \nu}$ is the outward unit normal to the
boundary with respect to $g$ (the time varying family of metrics), we have

\begin{prop}

For the unnormalized Ricci flow, the scalar curvature $R$ satisfies the following evolution equation.
\[
\left\{
\begin{array}{l}
\frac{\partial R}{\partial t}=\Delta R + R^2 \quad \mbox{in}\quad M\times\left(0,T\right),\\
\frac{\partial R}{\partial \nu}=k_g R-2k_g' \quad \mbox{on} \quad \partial M\times \left(0,T\right).
\end{array}
\right.
\]
Here the $'$ represents differentiation with respect to $t$.
\end{prop}
\bfseries \textit{Proof. } \normalfont A straightforward computation.

\hfill $\Box$

Since for problem (\ref{BVP0}) $k_g'=0$, as a consequence of the Maximum Principle we obtain
\begin{prop}
In problem (\ref{BVP2}), if $R\geq 0$ at $t=0$, it remains so.
\end{prop}
\section{Unnormalized flow: blow up of the curvature.}
\label{BlowUp}

The main technique we will use to study the limiting behavior of the Ricci Flow
is the analysis of the blow up limits. Therefore it is important to identify  
situations where the curvature blows up 
(in the case of the unnormalized flow). This is the content of the next proposition.

\begin{prop}
Assume $\chi(M)=1$ and $R>0$ through the flow. Assume that $k_g\geq 0$ and $k_g'\leq 0$.
Then the scalar curvature blows up in finite time.
\end{prop}
\bfseries \textit{Proof. }\normalfont 
Apply the Maximum Principle to the equation for the evolution of the curvature.

\hfill $\Box$

It can also be shown that

\begin{prop}
Assume $\chi(M)=1$ and that the geodesic curvature is nonpositive. Then the scalar curvature 
blows up in finite time.
\end{prop}

\bfseries \textit{Proof. }\normalfont 
Let $A(t)$ be the area of the manifold at time $t$. Then a computation shows that
\[
A'(t)=-\int_M R\, dV .
\]

This 
together with the fact $\chi(M)=1$ shows that
\[
A'(t)\leq -\epsilon, \quad \mbox{for some} \quad \epsilon>0
\]
Hence, the curvature cannot remain bounded after time $T=\frac{A(0)}{\epsilon}$ (Otherwise, by the
derivative estimates of Section \ref{DerivativeEstimates}, we would be able to continue the solution past $T$). 

\hfill $\Box$


\section{Normalized flow: Existence for all time.}\label{ExistenceForAllTime}
As we have seen, the equation 
we are studying has a solution for a short time. In general, we do not
expect this equation to have a solution for all time. 
In this section, under certain restrictions, we show longtime existence for the normalized flow
(in the rotationally symmetric case). To prove our result, we follow the ideas in \cite{Ham88}. 
Let us remark that the results in this section do not depend on derivative estimates on 
the curvature. 

We start with a definition.
\begin{defn}
We define the potential function as the solution to the problem
\[
\left\{
\begin{array}{l}
\Delta f= R-r\\
\frac{\partial f}{\partial \nu}=0 ,
\end{array}
\right.
\]
with mean value 0, where $r$ is the average of the scalar curvature, i.e., 
\[
r=\frac{\int_M R}{A\left(M\right)} \quad ,
\]
where $A\left(M\right)$ represents the area of $M$.
\end{defn}
We compute the evolution equation satisfied by the potential function.
\begin{lem}
The potential function satisfies an evolution equation
\[
\left\{
\begin{array}{l}
\frac{\partial f}{\partial t}=\Delta f+rf+\psi \quad \mbox{in}\quad M\times \left(0,T\right)\\
\frac{\partial f}{\partial \nu}=0 \quad \mbox{on} \quad \partial M \times \left(0,T\right) ,
\end{array}
\right.
\]
where $\psi$ satisfies
\[
\left\{
\begin{array}{l}
\Delta \psi=-r'\\
\frac{\partial \psi}{\partial \nu}=-\frac{\partial R}{\partial \nu}
=k\left(r-R\right)+2k' ,
\end{array}
\right.
\]
where $'$ denotes derivative with respect to time.
\end{lem}
\bfseries \textit{Proof.} \normalfont We compute
\[
\begin{array}{rcl}
\left(R-r\right)\Delta f +\Delta\left(\frac{\partial f}{\partial t}\right)
&=&\frac{\partial}{\partial t}\left(\Delta f\right)\\
&=&\frac{\partial}{\partial t}\left(R-r\right)\\
&=&\Delta R + R\left(R-r\right)-r'\\
&=& \Delta\left(\Delta f\right)+R\Delta f -r'.
\end{array}
\]
From this follows that
\[
\Delta\left(\frac{\partial f}{\partial t}-\Delta f-rf\right)=r' .
\]

On the other hand, if we write $\psi:=\frac{\partial f}{\partial t}-\Delta f-rf$,
we have,
\[
\frac{\partial \psi}{\partial \nu}=-\frac{\partial R}{\partial \nu}.
\]

\hfill $\Box$

Consider the function
\[
h=\Delta f+\left|\nabla f\right|^2 .
\]
Then we have that
\begin{lem}
$h$ satisfies an evolution equation
\[
\left\{
\begin{array}{l}
\frac{\partial h}{\partial t}=\Delta h -2\left|M_{ij}\right|^2+rh-r'-2\left<\nabla \psi,\nabla f\right>\\
h=R-r \quad \mbox{on}\quad \partial M\times \left(0,T\right) ,
\end{array}
\right.
\]
where $M_{ij}=\nabla_i\nabla_j f -\frac{1}{2}\Delta f g_{ij}$.
\end{lem}
\bfseries \textit{Proof.} \normalfont We compute
\[
\begin{array}{rcl}
\frac{\partial}{\partial t}\left|\nabla f\right|^2 &=&
2g^{ij}\left(\frac{\partial}{\partial t}\partial_i f\right)\partial_j f
+\left(\frac{\partial}{\partial t}g^{ij}\right)\partial_i f\partial_j f\\
&=& 2g^{ij}\partial_i\left(\Delta f+rf+\psi\right)\partial_j f+\left(R-r\right)g^{ij}\partial_i f\partial_j f\\
&=& 2\left<\nabla\Delta f,\nabla f\right>+2r\left|\nabla f\right|^2+\left<\nabla\psi,\nabla f\right>
+\left(R-r\right)\left|\nabla f\right|^2 .
\end{array}
\]

From the Weitzenb\"ock formula we obtain
\[
\begin{array}{rcl}
\frac{\partial}{\partial t}\left|\nabla f\right|^2
&=&\Delta\left|\nabla f\right|^2 -2\left|Hess f\right|^2-R\left|\nabla f\right|^2\\
&&
+2r\left|\nabla f\right|^2+\left<\nabla \psi,\nabla f\right>+\left(R-r\right)\left|\nabla f\right|^2 .
\end{array}
\]

To see that $h|_{\partial M}=R-r$ recall that $\left|\nabla f\right|^2=\left|\frac{\partial f}{\partial \nu}\right|^2=0$

\hfill $\Box$

From the formula in the previous Lemma we can deduce the following inequality:
\[
\begin{array}{rcl}
\frac{\partial h}{\partial t}&\leq&
\Delta h + rh -r'+2\left\|\nabla \psi\right\|_{\infty}\left|\nabla f\right|\\
&\leq& 
\Delta h +rh -r'+2\left\|\nabla \psi\right\|_{\infty}
\left(\left|\nabla f\right|+R-r+r\right)\\
&&\mbox{(since}\quad R>0\,\mbox{)}\\
&=&\Delta h +\left(r+2\left\|\nabla \psi\right\|_{\infty}\right)h-r'+2\left(r+\frac{1}{4}\right)\left\|\nabla \psi\right\|_{\infty} .
\end{array}
\]

To proceed we will need the following fact.
\begin{lem}
\label{Lenght}
Let $g_k=dr^2+f_k\left(r\right)^2d\omega$ be a sequence metrics on $B^2$. Let
$d_k=diam_{g_k}\left(M_k\right)$ and assume that there is a constant $C>0$ independent of $k$ such that
$d_k<C$, and there exist $\epsilon>0$ such that $k_{g_k}>-\epsilon$. Then
if $l_k(\partial B^2)\rightarrow 0$ (the length of the boundary
in the metric $g_k$) we must have $Vol_k\left(B^2\right)=:A_k\rightarrow 0$
\end{lem}
\bfseries \textit{Proof. }\normalfont Denote by $r_k=\frac{d_k}{2}$. If $r_k \rightarrow 0$, we are done. If not, let
$g_k$ be a subsequence such that $\lim_{k\rightarrow \infty} r_k=\alpha$. The hypothesis imply 
that $\alpha<\infty$.
Then, for any $0<\tau<\alpha$, we have that
\[
\frac{f_k'}{f_k}=k_{g_k} ,
\]
and hence,
\[
f_k(\alpha)=f_k(\tau)\exp\left(\int_{\tau}^{\alpha}k_{g_k}d\rho\right) .
\]

Notice that $k_{g_k}$ is decreasing, and then,
\[
\begin{array}{rcl}
f_k(\alpha) &\geq& f_k(\tau)\exp\left(\int_{\tau}^{\alpha}k_{g_k}(\alpha)\,d\rho\right)\\
&\geq& f_k(\tau)\exp\left(-\int_{\tau}^{\alpha}\epsilon\,d\rho \right)\\
&=& f_k(\tau)\exp\left[-\left(\alpha-\tau\right)\epsilon\right]\\
&\geq& f_k(\tau)\exp(-\epsilon\alpha).
\end{array}
\]

Since $l_k\left(B^2\right)\rightarrow 0$, $f_k\left(\alpha\right)\rightarrow 0$ also, and we have that
 $f_k(\tau)\rightarrow 0$ as $k\rightarrow \infty$. 
Therefore, for any $\eta>0$ we can find a $\tau_0$ such that $f_k(\tau_0)\leq \eta$,
and then choose $k$ large enough so that for any $\tau>\tau_0$ $f_k(\tau)<\eta$. This shows that,
\[
\begin{array}{rcl}
Vol(M_k)&\leq& 2\pi \int_0^{\alpha}f_k(\rho) \,d\rho\\
&\leq& 2\pi\eta\alpha ,
\end{array}
\]
which proves the lemma.

\hfill $\Box$

We are ready to prove now,
\begin{thm}
\label{longtime}
Assume $R\geq 0$ at $t=0$. Then solution to the normalized flow exists for all time.
\end{thm}
\bfseries\textit{Proof. }\normalfont The proof is going to be developed in two steps, and is
by contradiction. Assume the curvature blows up at some time $T<\infty$
\newline
\bfseries Step 1. \normalfont  \textit{We have that
\[
\int_0^T\left\|\nabla \psi\right\|_{\infty}\,dt<\infty
\]
}
First we notice that the function $|\nabla \psi|^2$ is subharmonic. Indeed, by the Weitzenb\"ock formula
\[
\begin{array}{rcl}
\Delta\left(\left|\nabla \psi\right|^2\right)&=&2\left|Hess\, \psi\right|^2+2\left<\nabla\left(\Delta \psi\right),\nabla \psi\right>
+R\left|\nabla \psi\right|^2\\
&\geq& 0 \quad \mbox{(notice that $\Delta \psi$ is constant in space)},
\end{array}
\]
so it attains its maximum at the boundary. Therefore,
\[
\begin{array}{rcl}
\left\|\nabla \psi\right\|&\leq& \left\|k_g\left(R-r\right)\right\|_{\infty,\partial M}+\left\|2k_g'\right\|_{\infty,\partial M}\\
&\leq& C\left\|R-r\right\|_{\infty,\partial M}+C 
\end{array}
\]
(of course we must show that $k_g'$ is bounded, see Remark \ref{geodesiccurvature} at the end of this section).

Hence, if $\left\|\nabla \psi\right\|_{\infty}$ is not integrable in $\left[0,T\right]$, neither is $R|_{\partial M}$. In such a case,
since the conformal factor is given by
\[
u\left(x,t\right)=\exp\left[\int_0^t r\left(\tau\right)-R\left(x,\tau\right)\,d\tau\right] ,
\]
$l\left(\partial M\right)\rightarrow 0$ as $t\rightarrow T$. But $R>0$ and $k_g\geq0$, and hence by Lemma \ref{Lenght}, we must
have $A\left(M\right)\rightarrow 0$, which contradicts the fact that the normalized flow keeps the volume constant.
\newline
\bfseries Step 2. \normalfont \textit{If the normalized flow does not exist for all time, we have that for the unnormalized flow
the curvature has maximum blow up rate in the boundary. 
}

As we already showed, the following inequality holds
\[
\frac{\partial h}{\partial t}\leq \Delta h+\left(r+2\left\|\nabla \psi
\right\|_{\infty}\right)h-r'+2\left(r+\frac{1}{4}\right)\left\|\nabla \psi\right\|_{\infty} .
\]

Let $c(t)=\int_0^t \left(r+2\left\|\nabla \psi\right\|_{\infty}\right)\,dt$, and define
\[
w=\exp\left(-c(t)\right)h ,
\]
then
\[ 
\frac{\partial w}{\partial t}\leq \Delta w -\exp\left(-c(t)\right)\left[r'-2\left\|\nabla \psi\right\|_{\infty}\left(r+\frac{1}{4}\right)\right] .
\]

Using the fact that $c(t)$ is bounded, the inequalities
\[
\begin{array}{rcl}
\left|r'\right|&\leq&\int_{\partial M} \left|\frac{\partial R}{\partial \nu}\right| \,d\sigma \\
&\leq& C\left[k_g\left\|R-r\right\|_{\infty,\partial M}+2\left\|k_g'\right\|_{\infty,\partial M}\right]\\
&\leq& C\left\|R-r\right\|_{\infty,\partial M}+C , 
\end{array}
\]
and
\[
\left|\nabla \psi\right|\leq C\left|R-r\right|_{\partial M}+C,
\]
the Maximum Principle shows that 
\[
h(t)\leq C\hat{R}_{max}(t)+C,
\]
where
\[
\hat{R}_{max}(t)=\max_{\left(x,\tau\right)\in\partial M\times\left[0,t\right]}R(x,\tau) .
\]

This proves the claim since $h\left(t\right)\geq R\left(t,x\right)-r\left(t\right)$, and 
$r\left(t\right)$ is bounded.

The last claim shows that $\int_0^T R_{\max}(\tau)\,d\tau<\infty$, where
$[0,T)$ is the maximum interval of existence for the Ricci Flow. Again,
we use the fact that the conformal factor $u$ satisfies the identity,
\[
u\left(x,t\right)=\exp\left[\int_0^t r\left(\tau\right)-R\left(x,\tau\right)\,d\tau\right] .
\]

It follows that there is a constant $\delta>0$ such that
$u\geq \delta$, and hence we can continue the solution past $T$ (since $u$
is bounded, by Theorem 1.3 in chapter III of \cite{DiB}, $u$ is H\"older
continuous. The rest follows from a standard bootstrapping argument).

\hfill $\Box$

\begin{rem}
\label{geodesiccurvature}
Assume that the normalized flow becomes singular at time $T_0<\infty$ we will show that $k_g'$ 
remains bounded. To see this, consider the unnormalized flow. Then $T_0$ corresponds to
the blow up time $T<\infty$ for the unnormalized flow. If $A\left(t\right)$
represents the area of the surface at time $t$, we cannot have $A\left(T\right)=0$ (because,
by using Gauss-Bonnet and the fact that we keep $k_g=k_0$ constant for the unnormalized flow, 
it can be shown that $A\left(t\right)\leq C\left(T-t\right)$, and this would imply that the
normalized flow exists for all time). Therefore, there is $\epsilon>0$ such that $A\left(t\right)\geq \epsilon$
for $0\leq t <T$.

On the other hand, if $\tilde{t}$ and $t$ represent the time parameter for the normalized and
unnormalized flow respectively, and $\phi = 1/A$, we have
\[
\frac{d}{d\tilde{t}}k_{g}=\frac{d}{dt}\left(\frac{k_0}{\sqrt{\phi\left(t\right)}}\right)\slash\frac{d\tilde{t}}{dt}
=-\frac{k_0}{2\phi^{\frac{5}{2}}}\phi' .
\]

Hence all we must show is that $\phi'=\left(\frac{1}{A}\right)'$ is bounded. But this is rather easy to see,
since
\[
A'\left(t\right)=-\int_{M}R\,dV =-2\pi +\int_{\partial M} k_0\,d\sigma , 
\]
which is clearly bounded.

\end{rem}
\hfill $\Box$

Let us observe that the long time existence result conveys some interesting information,
namely
\begin{prop}
For the normalized flow, $k_g\rightarrow 0$ exponentially.
\end{prop}

\section{Unnormalized flow: Derivative estimates.}\label{DerivativeEstimates}

To produce estimates on the derivatives of the curvature from bounds on the curvature
we will fix a collar of the boundary at time $t=0$. When working with a rotationally symmetric
metric, if we have an upper bound on the curvature and a bound on the geodesic curvature
of the boundary, the Laplacian Comparison Theorem gives an estimate on the radius of the ball with
the endowed metric. This fact allows us to choose a rotationally symmetric 
collar $V=\partial M \times [0,c)$ such that $\rho_{t}\left(\cdot, \partial M\right)$, the
distance to the boundary at time $t$, is smooth, and the distance at time $t$ from the
boundary $\partial M=\partial M\times \left\{0\right\}$ and the inner boundary
of the collar $\partial M\times \left\{c\right\}$ is bounded below uniformly by $\frac{c}{2}$.

Also, by the interior derivative estimates of Shi, on the inner boundary of the collar for $t>0$ we 
can assume bounds on all the derivatives of the curvature in terms of a bound of the curvature. 
All these said, we set ourselves to the task of estimating the derivatives $\frac{\partial^n R}{\partial \nu^n}$,
where by $\frac{\partial}{\partial \nu}$ we denote the radial unit vector which is well defined 
in $V$ at all times.  

\subsection{First Order Estimates.}\label{FirstOrderDerivatives}
The following Theorem and its proof are along the lines of Theorem 7.1 in \cite{Ham95} .
\begin{thm}
There exist constants $C_1$ for $R\geq 1$ such that if the curvature is bounded
\[
|R|\leq M ,
\]
up to time $t$ with $0<t\leq\frac{1}{M}$ then the covariant derivative of the curvature is bounded
\[
|\frac{\partial R}{\partial \nu}|\leq C_1\frac{M}{\sqrt{t}} .
\]
\end{thm}
\bfseries \textit{Proof.} \normalfont 
At the boundary we have the identity
\[
\frac{\partial R}{\partial\nu}=k_g R .
\]

Define the quantity $F=t|\frac{\partial R}{\partial\nu}|^2+AR^2$, where $A$ is a constant to be chosen.
In the interior we have 
\[
\begin{array}{rcl}
\frac{\partial F}{\partial t}&=& (R_{\nu})^2+2t(R_{\nu})_tR_{\nu}+2ARR_t\\
&=& (R_{\nu})^2+2tR_{\nu}\left[RR_{\nu}+(R_t)_{\nu}\right]+2AR\left[\Delta R +R^2\right]\\
&=& (R_{\nu})^2+2tR_{\nu}\left[(\Delta R)_{\nu}+2RR_{\nu}+RR_{\nu}\right]+2AR\left[\Delta R +R^2\right]\\
&\leq&\Delta \left(t(R_{\nu})^2+AR^2\right)-2\left(t|\nabla R_{\nu}|^2+A|\nabla R|^2\right)
+CtR(R_{\nu})^2+2AR^3\\
&\leq&\Delta F +\left(CtR-2(A-1)\right)|\nabla R|^2+2AR^3 .
\end{array}
\]

By choosing $A$ big enough we get the inequality,
\[
\frac{\partial F}{\partial t}\leq \Delta F+2AR^3 .
\]

At the boundary of $M$ we have (and we may assume a similar estimate on the inner boundary
of the collar we fixed from the beginning),
\[
F=tk_g^2 R^2+A R^2 \leq \tilde{C}tM^3+\tilde{C}M^2 ,
\]
for some constant $\tilde{C}$. Then the maximum principle, as long as $tM\leq 1$, yields 
\[
t|\frac{\partial R}{\partial \nu}|^2\leq F\leq CM^2 ,
\]
for some constant $C_1$.

\hfill $\Box$

\subsection{Second Order Estimates.}\label{SecondOrderEstimates}
We show how to estimate $\frac{\partial^2 R}{\partial \nu^2}$. Higher derivatives
follow a similar procedure. 

Let 
\[
\rho(P,t)=d_t(P,\partial M)
\]
be the distance from $P$ to the boundary at time $t$ (notice we have dropped the index that indicates 
the time dependence of $\rho$, but it is always well understood that this is so),  and define
\[
F=\exp(k_0\rho)R \quad \mbox{where $k_0$ is the geodesic curvature of $\partial M$}
\]
(we want to remark that if $k_0=0$ all these preparations become unnecessary).

A simple computation shows that $F_{\nu}=0$ at the boundary. $F$ satisfies the evolution equation
\begin{eqnarray*}
\lefteqn{F_t=\Delta F-2\left(\frac{\nabla\exp(k_0\rho)}{\exp(k_0\rho)}\right)\cdot \nabla F}\\
&&+2\left(\frac{|\nabla\exp(k_0\rho)|^2-\Delta\exp(k_0\rho)}{\exp(k_0\rho)}+\frac{k_0\rho'}{2}\right)F+\frac{F^2}{\exp(k_0\rho)},
\end{eqnarray*}
where $'$ denotes differentiation with respect to $t$. To simplify, we rewrite this last expression as
\[
\left\{
\begin{array}{l}
F_t=\Delta F+ B\cdot \nabla F +CF +EF^2 \quad \mbox{in} \quad B^2\times\left(0,T\right)\\
F_{\nu}=0 \quad \mbox{on}\quad \partial B^2\times\left(0,T\right), 
\end{array}
\right.
\]

Where $B,C,E$ have the obvious meanings. Differentiating with respect to $\nu$ and writing $w=F_{\nu}$ we get,
\begin{eqnarray*}
Lw\equiv\lefteqn{w_t-\Delta w -B\cdot \nabla w}\\
&&=RF_{\nu}+B_{\nu}\cdot \nabla F-C_{\nu}F+E_{\nu}F^2-(C+EF)w .
\end{eqnarray*}

Write
\[
\mathcal{G}=RF_{\nu}+B_{\nu}\cdot \nabla F-C_{\nu}F+E_{\nu}F^2-(C+EF)w , 
\]
and the boundary condition is $w=0$

To estimate derivatives of $w$, we want to use the theory of Fundamental Solutions for linear parabolic
equations. In order to do so, we are going to compute the coefficients of the operator $L$ in Fermi
coordinates.

Pick a point $P\in \partial M$. Assume that $P$
admits Fermi coordinates $(\mathcal{F}_{\epsilon,\delta},\phi)$ in a neighborhood
\[
\mathcal{F}_{\epsilon,\delta}=\left\{(x,s)\in \mathbf{R}^2:\,-\delta<x<\delta \quad 0\leq s<\epsilon  \right\}
\]
i.e., the $s$-coordinate represents the distance to the boundary, and $x$ is the coordinate on the curves ``parallel''
to the boundary. By $k(\sigma)$ we will denote the geodesic curvature of
\[
\partial M_{\sigma}=\left\{P\in M\,:\, \pi_2\left(\phi^{-1}(P)\right)=\sigma \right\} , 
\]
where $\pi_2\,:\,\mathbf{R}^2\rightarrow \mathbf{R}$ is the projection onto the second coordinate.

In Fermi coordinates the metric can be described explicitly in terms of the relevant geometric quantities. This is the
content of the next proposition.

\begin{prop}
\label{coeff}
In Fermi coordinates,
\newline
(1) $g=ds^2+\exp\left(-2\int_0^s k(\sigma)\,d\sigma\right)dx^2$;
\newline
(2) $\frac{\partial k}{\partial s}=k^2+\frac{R}{2}$;
\newline
(3) $\sqrt{g}=\exp\left(-\int_0^s k(\sigma)\,d\sigma\right)$;
\newline
where $s$ is the distance to the boundary.
\end{prop}
\bfseries \textit{Proof. } \normalfont Let us prove (1). Just notice that,
\[
\begin{array}{rcl}
\frac{\partial}{\partial s}g\left(\frac{\partial}{\partial x},\frac{\partial}{\partial x}\right)
&=&2g\left(\nabla_{\frac{\partial}{\partial s}}\frac{\partial}{\partial x},\frac{\partial}{\partial x}\right)\\
&=&-2g\left(\nabla_{\frac{\partial}{\partial x}}\frac{\partial}{\partial s},\frac{\partial}{\partial x}\right)=-2k(s)
g\left(\frac{\partial}{\partial x},\frac{\partial}{\partial x}\right) .
\end{array}
\]
The result follows by integration and recalling that at $s=0$
$g\left(\frac{\partial}{\partial x},\frac{\partial}{\partial x}\right)=1$. We go for (2). On one hand we have,
\[
k(s)g\left(\frac{\partial}{\partial x},\frac{\partial}{\partial x}\right)=
-g\left(\nabla_{\frac{\partial}{\partial x}}\frac{\partial}{\partial s},\frac{\partial}{\partial x}\right) .
\]

Differentiate the previous equation with respect to $s$. First the left hand side,
\[
\begin{array}{rcl}
\frac{\partial}{\partial s}\left[k(s)g\left(\frac{\partial}{\partial x},\frac{\partial}{\partial x}\right)\right]
&=& k'(s)g\left(\frac{\partial}{\partial x},\frac{\partial}{\partial x}\right)+
2k(s)g\left(\nabla_{\frac{\partial}{\partial s}}\frac{\partial}{\partial x},\frac{\partial}{\partial x}\right)\\
&=& k'(s)g\left(\frac{\partial}{\partial x},\frac{\partial}{\partial x}\right)
-2[k(s)]^2g\left(\frac{\partial}{\partial x},\frac{\partial}{\partial x}\right) .
\end{array}
\] 

For the right hand side we have
\[
\begin{array}{rcl}
\frac{\partial}{\partial s}g\left(\nabla_{\frac{\partial}{\partial x}}\frac{\partial}{\partial s},\frac{\partial}{\partial x}\right)
&=& g\left(\nabla_{\frac{\partial}{\partial s}}\nabla_{\frac{\partial}{\partial x}}\frac{\partial}{\partial x}\right)
+g\left(\nabla_{\frac{\partial}{\partial x}}\frac{\partial}{\partial s},\nabla_{\frac{\partial}{\partial s}}\frac{\partial}{\partial x}\right)\\
&=& -\frac{R}{2}g\left(\frac{\partial}{\partial x},\frac{\partial}{\partial x}\right)+
[k(s)]^2g\left(\frac{\partial}{\partial x},\frac{\partial}{\partial x}\right) .
\end{array}
\]

Combining the two previous computations give the result. (3) follows from (1).

\hfill $\Box$

As we said before, 
we are going to study the operator $L$ in Fermi coordinates. First we study the principal part which corresponds to the Laplacian.
 Recall that the Laplacian in local coordinates is given by
\[
\Delta_g = \frac{1}{\sqrt{g}}\partial_i\left(\sqrt{g}g^{ij}\partial_j\right) .
\]

Choose Fermi coordinates at time $t=0$. In this coordinates the Laplacian is given by
\begin{eqnarray*}
\lefteqn{\Delta =
\frac{\partial^2}{\partial^2 s}+\exp\left(2\int_0^s k(\sigma)\,d\sigma\right)\frac{\partial^2}{\partial^2 x}}\\
&&+\exp\left(\int_0^s k(\sigma)\,d\sigma\right)
\frac{\partial}{\partial s}\left[\exp\left(-\int_0^s k(\sigma)\,d\sigma\right)\right]\frac{\partial}{\partial s} .
\end{eqnarray*}

Hence, since $\Delta_{g(t)}=\exp\left(\int_0^t R(s,\tau)\,d\tau\right)\Delta$
(that the conformal factor is given by $u=\exp\left(-\int_0^t R\left(s,\tau\right)\,d\tau\right)$ and 
in dimension 2, if $h=\lambda g$ then
$\Delta_h=\lambda^{-1}\Delta_g$), the functions whose H\"older norms
we have to estimate are
\newline
(1) $\alpha(s,t)=\exp\left(\int_0^t R(s,\tau)\,d\tau\right)$.
\newline
(2) $\beta(s,t)=\exp\left(\int_0^t R(s,\tau)\,d\tau\right)\exp\left(2\int_0^s k(\sigma)\,d\sigma\right)$.
\newline
(3) $\gamma(s,t)=k(s)\exp\left(\int_0^t R(s,\tau)\,d\tau\right)$.

\medskip
We will check that these functions are differentiable and that their derivatives can be bounded in terms of
bounds on the curvature and its first derivative (which we already know how to bound).
 We do it for (3), since it is the same for the other functions. Notice that by the rotational
symmetry of the metric, which is preserved by the flow, the $\frac{\partial}{\partial x}$ derivatives
are 0. So we only have to compute the $t$ and $s$-derivatives. Let us compute first the $s$-derivatives:
\[
\begin{array}{rcl}
\gamma_s&=&k'(s)\exp\left(\int_0^tR(s,\tau)\,d\tau\right)+\\
&& k(s)\left(\int_0^t R_s(s,\tau)\,d\tau\right)\exp\left(\int_0^tR(s,\tau)\,d\tau\right) .
\end{array}
\]

By taking into account (2) from proposition
 (\ref{coeff}), we have an estimate on $\gamma_s$ in terms of the curvature, its first derivative, and
the geodesic curvature. For the time derivative we have,
\[
\gamma_t = R(s,t)k(s)\exp\left(\int_0^t R(s,\tau)\,d\tau\right),
\]
and hence the same conclusion.

Now we work with the lower order terms. We have the following list of explicit formulae
\begin{prop} The following formulae hold,
\newline
(1) $|\nabla_{g(t)} \exp(k_0\rho)|^2=k_0^2\exp\left(2k_0\rho\right)$,
\newline
(2) $\Delta_{g(t)}\exp(k_0\rho)=\exp\left(\int_0^t R(\sigma,\tau)\,d\tau\right)\left[k_0^2\exp\left(k_0\rho\right)+k_0k(s)\exp\left(k_0\rho\right)\right]$,
\newline
(3) $\rho(s,t)=\int_0^s\exp\left(-\int_0^t R(\sigma,\tau)\,d\tau\right)\,d\sigma$.
\end{prop}
With the aid of the previous  
 proposition, we can compute the terms $B,C,E$, and show that they satisfy a Lipschitz condition, and that their Lipschitz constants
depend only on bounds on the curvature. Also, it can be easily checked that $B_{\nu},C_{\nu}, E_{\nu}$ are bounded continous functions and their
bounds depend only on bounds on the curvature and its first derivative which we can bound
in terms of the curvature.

To take our next step we need the following proposition.
\begin{prop}
\label{interiorestimates}
In the strip $\left[-\delta,\delta\right]\times\left[\frac{3}{4}\epsilon, \epsilon\right]$ we have uniform bounds on all the
covariant derivatives of the curvature, and these bounds depend only on $\epsilon$, $t$ and a bound on the curvature
\end{prop}
\bfseries \textit{Proof. } \normalfont The proof follows from the local interior estimates (See Theorem 13.1 in \cite{Ham95}).

\hfill $\Box$

Fix a function $\psi:\mathbf{R}_+\longrightarrow [0,1]$, such that
\[
\psi(u) =\left\{
\begin{array}{l}
1 \quad \mbox{if} \quad u\in [0,1]\\
0 \quad \mbox{if}  \quad u\in[2,\infty] \quad .
\end{array}
\right.
\]

Define,
\[
\chi(s)=\psi\left[\frac{8\left(\epsilon-s\right)}{\epsilon}\right] ,
\]

Notice that there is a constant $A$ such that 
\[
\left|\nabla \chi\right|\leq\frac{A}{\epsilon} \quad \mbox{and} \quad \left|\nabla^2 \chi\right|\leq \frac{A^2}{\epsilon^2} .
\]

 Let
\[
v=w-\chi w .
\]

A simple computation shows that
\begin{eqnarray*}
\lefteqn{v_t-\Delta v-B\cdot\nabla v}\\
&&=\mathcal{G}-\chi\mathcal{G}+w\Delta \chi+2\nabla\chi\cdot \nabla w +Bw\nabla \chi-Gv .
\end{eqnarray*}

From now on we will use the following notation
\[
\mathcal{W}=\mathcal{G}-\chi\mathcal{G}+w\Delta \chi+2\nabla\chi\cdot \nabla w +Bw\nabla \chi-Gv .
\]

Since for $s\leq\frac{3}{4}\epsilon$ we have $\chi\equiv 0$, the left hand side of the previous
equation (by Proposition \ref{interiorestimates}) 
is bounded in terms of bounds on the curvature and $\epsilon$. This shows that
the left hand side is a bounded continous function on the segment $\{0\leq s<\epsilon\}$. Pick 
any $\theta>0$, and set $f=v|_{t=\theta}$.
Let $\Gamma(s,t;\zeta,\tau)$ be the Green's function of the operator $L$ in the layer
$\left\{0<s<\epsilon\right\}\times [\theta,T]$
Since $v=0$ at $s=0$ and $s=\epsilon$, we can represent $v$ as follows
\begin{equation}
\label{representation}
v(s,t)=\int_{\theta}^t\int_0^{\epsilon} \Gamma(s,\zeta;t,\tau)\mathcal{W}(\zeta,\tau)\,d\zeta\,d\tau
+\int_0^{\epsilon} \Gamma(s,\zeta;t,\theta)f(\zeta)\,d\zeta .
\end{equation}

It is well known that we can differentiate $v$ with respect to $s$, and the 
derivative is given by direct differentiation under the integral sign. That is to say,
\begin{eqnarray*}
\lefteqn{\frac{\partial v}{\partial s}=
\int_{\theta}^t\int_0^{\epsilon} \frac{\partial}{\partial s}\Gamma(s,t;\zeta,\tau)\mathcal{W}(\zeta,\tau)\,d\zeta\,d\tau}\\
&&\qquad\qquad\qquad+\int_0^{\epsilon}\frac{\partial}{\partial s}\Gamma(s,t;\zeta,\tau)f(\zeta)\,d\zeta .
\end{eqnarray*}

We are ready to produce bounds on $\frac{\partial v}{\partial s}$. First of all, it is well known that
\[
\begin{array}{c}
\Gamma\left(s,t;\xi,\tau\right)\leq  
C\left(t-\tau\right)^{-\frac{1}{2}}\exp\left[-C\frac{\left|s-\xi\right|^2}{t-\tau}\right]\\
\frac{\partial}{\partial s}\Gamma\left(s,t;\xi,\tau\right)
\leq C\left(t-\tau\right)^{-\frac{1+1}{2}}\exp\left[-C\frac{\left|s-\xi\right|^2}{t-\tau}\right] ,
\end{array}
\]
where $C$ depends only on the Lipschitz bounds on the coefficients of $L$.
Therefore, if 
\[
\left|f\left(\xi\right)\right|\leq M_0
\quad \mbox{and} \quad 
\left|\mathcal{W}\left(\xi,\tau\right)\right|\leq M_1 , 
\]
writing
\[
\int_{\theta}^\infty\int_0^{\infty} 
 C\left(t-\tau\right)^{-\frac{1+1}{2}}\exp\left[-C\frac{\left|s-\xi\right|^2}{t-\tau}\right] \,d\xi\,d\tau
=K_1
\]
and
\[
\int_0^{\infty}  C\left(t-\tau\right)^{-\frac{1}{2}}
\exp\left[-C\frac{\left|s-\xi\right|^2}{t-\tau}\right] \,d\xi
=K_0 ,
\]
we find that
\[
\left|\frac{\partial v}{\partial s}\right|\leq K_1 M_1 + \frac{M_0 K_0}{\left(t-\theta\right)^{\frac{1}{2}}} .
\]

Finally, notice that $M_0, M_1, K_0, K_1$ depend only on bounds on the curvature and its first derivative. Since
\[
\begin{array}{rcl}
\frac{\partial v}{\partial s}\left(x,t\right)&=&
\frac{\partial^2 w}{\partial s\partial\nu}\left(x,t\right)
\quad\left(\mbox{recall that}\, 
\frac{\partial}{\partial \nu}= 
g\left(\frac{\partial}{\partial s},\frac{\partial}{\partial s}\right)^{-1}\frac{\partial}{\partial s}\right)\\
&=&\exp\left(-\frac{1}{2}\int_0^t R\left(x,\tau\right)\,d\tau\right)\frac{\partial^2 w}{\partial \nu^2}\left(x,t\right)\\
&=&\exp\left(-\frac{1}{2}\int_0^t R\left(x,\tau\right)\,d\tau\right)
\left[\exp\left(k_0\rho\right)\frac{\partial^2 R}{\partial \nu^2}+\left(k_0^2-k_0\right)\exp
\left(k_0\rho\right)\frac{\partial R}{\partial \nu}\right] .
\end{array}
\]

We have proved the following
\begin{thm} 
Second partial derivatives of $R$ (in a Fermi coordinate neighborhood)
can be bounded up to the boundary, and the bounds only depend on bounds on $R$, 
the geodesic curvature of the boundary, the time $t>0$ up to
where we are bounding, and $\epsilon>0$, the ``size'' of the Fermi coordinate neighborhood at $t=0$.
\end{thm}

\section{Subsequential Convergence: the case of convex boundary ($k_0\geq 0$).}
\label{Perelman}

Our purpose in this section is to show that when $k_0\geq 0$, the curvature is uniformized by the (unnormalized) flow in the
sense that if $[0,T)$ is the maximum interval of existence of the solution for the unnormalized
flow, then along a sequence of times $t_k\rightarrow T$ 
\[
\lim_{k\rightarrow \infty} \frac{R_{max}\left(t_k\right)}{R_{min}\left(t_k\right)}=1 ,
\]
where $R_{max}\left(t\right)=\max_{x\in M}R\left(x,t\right)$, and $R_{min}\left(t\right)$ is
defined analogously. 

To accomplish our goal we will use blow up analysis. First we describe how to form a blow up limit.
Let $\left(0,T\right)$ be the maximum time of existence of the Ricci flow. Then we can find a 
sequence of times $t_i\rightarrow T$ and points $p_i\in M$ such that
\[
\lambda_i:=R\left(p_i,t_i\right)=\max_{M\times \left[0,t_i\right]}R\left(x,t\right) ,
\]
and we have $\lambda_i\rightarrow \infty$. Define the dilations

\[
g_i\left(t\right):=\lambda_i g\left(t_i+\frac{t}{\lambda_i}\right) \quad 
-\lambda_i t_i\leq t<\lambda_i\left(T-t_i\right) .
\]

Then, using the derivative estimates (plus injectivity radius and Fermi inradius assumptions,
easily seen to be fulfilled in the rotationally symmetric case via the 
Laplacian Comparison Theorem), we can find a subsequence of times, which we
denote again by $t_i\rightarrow T$, such that $\left(M,g_i\left(t\right), p_i\right)$ converges smoothly
to a solution to the Ricci flow $\left(M_{\infty},g_{\infty}\left(t\right),p_{\infty}\right)$.
 
The following theorem is the key to our approach.
\begin{thm}[See Theorems 26.1 and 26.3 in \cite{Ham95}]
\label{limitingbehavior}
There are two possible blow up limits for solutions of \ref{BVP0} with $R>0$. If the blow up limit is compact, then it is 
the round hemisphere $S^2_+$ with totally geodesic boundary. If the blow up limit
is noncompact, then it (or its double) is isometric to the cigar metric in $\mathbf{R}^2$ 
\end{thm}

Recall that the cigar metric in $\mathbf{R}^2$ is given by
\begin{equation}
\label{cigar}
ds^2=\frac{dx^2+dy^2}{1+x^2+y^2}
\end{equation}

\bfseries \textit{Proof. } \normalfont The main observation here is that the geodesic curvature of a rescaled
metric $\lambda^2 g$ is given by
\[
k_{\lambda^2 g}=\frac{k_g}{\lambda} .
\]

Therefore the blow up limits (as long as the geodesic curvature is kept uniformly bounded throughout the flow) have totally
geodesic boundary.   
Once this has been established, it is not difficult to show (by following the corresponding arguments in the boundaryless case)
that the conclusions of the Theorem hold.

\hfill $\Box$

According to the previous Theorem, in order to show that the flow is uniformizing the curvature, all 
 we need to show is that among all the possibilities, the only possible blow up limit is $S^2_+$. 
Since we already know that the only other
possibility is the soliton metric in $\mathbf{R}^2$ (by Theorem \ref{limitingbehavior}),
all we have to do is to rule
out this case, and here is where Perelman's work (see \cite{Per1}) comes handy. We will follow his ideas
and will apply them in a less general framework.

Following Perelman, first we define the following functional,
\begin{equation}
\label{perelman0}
\mathcal{F}\left(g,f\right)=\int_M\left(R+\left|\nabla f\right|^2\right)\exp\left(-f\right)\,dV
\end{equation}
and compute its first variation.

\begin{prop}
Let $\delta g_{ij}=v_{ij}, \delta f = h, g^{ij}v_{ij}=v$. Then we have,
\begin{equation}
\label{firstvariation}
\begin{array}{rcl}
\delta \mathcal{F}
&=&\int_M\exp\left(-f\right)\left[-v_{ij}\left(R_{ij}+\nabla_i\nabla_j f\right)
+\right(\frac{v}{2}-h\left)\left(2\Delta f-\left|\nabla f\right|^2+R\right)\right]\\
&&-\int_{\partial M}\left(\frac{\partial v}{\partial \nu} +v\frac{\partial f}{\partial \nu}\right)\exp\left(-f\right)\,d\sigma+\\
&&\int_{\partial M}\exp\left(-f\right)\nabla_i v_{ij}\nu^j\,d\sigma
-\int_{\partial M} \nabla_j\exp\left(-f\right)v_{ij}\nu^i\,d\sigma .
\end{array}
\end{equation}
here, $\frac{\partial}{\partial \nu}=\left\{\nu^i\right\}$ is the outward unit normal to $\partial M$ with respect to $g$ (notice that then 
$\frac{\partial}{\partial\nu}$ depends on time), 
$\nabla$ represents
covariant differentiation
with respect to the metric $g$, and $v_{\partial}/2$ represents the  variation of the volume element of $\partial M$
induced by $v_{ij}$.
\end{prop}
\bfseries\textit{Proof. }\normalfont Equation (\ref{firstvariation}) follows from the following formulas
\begin{equation}
\label{variation1}
\begin{array}{c}
\delta R = -\Delta v+\nabla_i\nabla_j v_{ij}-R_{ij}v_{ij}\\
\delta \left|\nabla f\right|^2=-v^{ij}\nabla_i f\nabla_j f+\left<\nabla f,\nabla h \right>\\
\delta \left(\exp\left(-f\right)\,dV\right)=\left(\frac{v}{2}-h\right)\exp\left(-f\right)\,dV ,
\end{array} 
\end{equation}
and integration by parts ($\Delta$ is the Laplacian operator of the metric $g$).

\hfill $\Box$
 
Consider the evolution equations

\begin{equation}
\label{evolution1}
\left\{
\begin{array}{l}
g_t=-2\left(Ric\left(g\right)+Hess f\right)\quad \mbox{in} \quad M\times \left(0,T\right)\\
k_g =k_0 \quad \mbox{on} \quad \quad \partial M \times \left(0,T\right)\\
f_t=-R-\Delta f \quad \mbox{in} \quad M \times \left(0,T\right)\\
\frac{\partial f}{\partial \nu}=0 \quad \mbox{on} \quad \partial M\times \left(0,T\right) .
\end{array}
\right.
\end{equation}

We have the following computation.
\begin{prop}
\label{monotonicity}
\begin{equation}
\label{firstvariation2}
\begin{array}{rcl}
\mathcal{F}_t&=&2\int_M \left|R_{ij}+\nabla_i\nabla_j f\right|^2\exp\left(-f\right)\,dV \\
&&+2\int_{\partial M}k_0 R\exp\left(-f\right)\,dA+2\int_{\partial M}k_0 \left|\nabla^{\partial}f\right|^2\exp\left(-f\right)\,dA.
\end{array}
\end{equation}
\end{prop}
\bfseries\textit{Proof. }\normalfont
We compute each of the boundary integrals in (\ref{variation1}) using the evolution equations (\ref{evolution1})
\[
\int_{\partial M}\exp\left(-f\right)\nabla_i v_{ij}\nu^j\,dA =
\int_{\partial M}k_0 R \exp\left(-f\right)\,dA,
\]

\[
-\int_{\partial M}\nabla_j\exp\left(-f\right)v_{ij}\nu^i\,dA =2\int_{\partial M}\exp\left(-f\right)k_0\left|\nabla^{\partial}f\right|^2\,dA,
\]

\[
-\int_{\partial M}\nabla_{\nu} \exp\left(-f\right)v_{\nu\nu}\,dA= \int_{\partial M}\nabla_{\nu}v \exp\left(-f\right)\,dA=0.
\] 

The monotonicity formula then follows. 

\hfill $\Box$

We will need the following more sophisticated version of $\mathcal{F}$,
\[
\mathcal{W}\left(g,f,\tau\right)=\int_M 
\left[\tau\left(\left|\nabla f\right|^2+R\right)+f-2\right]
\left(4\pi\tau\right)^{-1}\exp\left(-f\right)\,dV,
\]
restricted to $f$ satisfying 
\[
\int_M \left(4\pi\tau\right)^{-\frac{1}{2}}\exp\left(-f\right)\,dV=1 .
\]

As pointed out by Kleiner and Lott, by making the substitution $\Phi=\exp\left(-\frac{f}{2}\right)$,
we get the functional
\begin{equation}
\label{perelman2}
\mathcal{W}\left(g,\Phi,\tau\right)=
\left(4\pi\tau\right)^{-1}\int_M 
\left[4\tau\left|\nabla \Phi\right|^2+\left(\tau R-2\log \Phi-2\right)\Phi^2\right]\,dV
\end{equation}
which has been extensively studied by Rothaus (See 
\cite{Rot1},\cite{Rot2}) in a domain with boundary $\Omega$
under the further restriction that the infimum is taken over smooth functions
vanishing in the boundary. Using his methods is not difficult to show that for our modified functional
there is a minimizer. 

Just as before, given a time dependent family of metrics evolving by the Ricci flow, i.e., satisfying
\[
\left\{
\begin{array}{l}
g_t=-Rg \quad \mbox{in} \quad M\times \left(0,T\right) \\
k_g=k_0 \quad \mbox{on}\quad \partial M\times\left(0,T\right)
\end{array}
\right.
\]
and  $f\left(\cdot,t\right)$ satysfying the backward heat equation 
\[
f_t=-\Delta f +\left|\nabla f\right|^2-R + \frac{1}{\tau} ,\quad \tau_t=-1 ,
\] 
plus 
the boundary condition
\[
\frac{\partial f}{\partial \nu}=0 .
\]

We have the following monotonicity formula. Its proof is 
an adaptation of the calculations in the proof of Proposition \ref{monotonicity}.

\begin{thm}
\label{monotonicity2}
\[
\begin{array}{rcl}
\frac{d\mathcal{W}}{dt}&=&\int_M 2\tau\left|R_{ij}+
\nabla_i\nabla_j f-\frac{1}{2\tau}g_{ij}\right|^2\left(4\pi\tau\right)^{-1}
\exp\left(-f\right)\,dV\\
&&+2\tau\left(\int_{\partial M}k_0 R\exp\left(-f\right)\,dA+\int_{\partial M}k_0 \left|\nabla^{\partial}f\right|^2\exp\left(-f\right)\,dA\right).
\end{array}
\]
\end{thm}

It follows immediately that the quantitity 
$\mu\left(g,\tau\right)=\inf_{f\in \mathcal{C}^{\infty}} \mathcal{W}\left(g_{ij},f,\tau\right)$ is
increasing in $\tau$ along the flow (of course under the assumption that $k_0\geq 0$).

\subsection{The Argument (Proof of Theorem \ref{convexboundarycase}). }

First we show an interesting and well known property of
the soliton metric, but before we start, let us set some notation

\bfseries Notation. \normalfont
The annulus in the manifold $\left(M,g\right)$ of inner radius $r_1>0$ and outer radius $r_2>r_1$ 
(in the metric $g$) will
be denoted by $A\left(r_1,r_2\right)$, i.e.,
\[
A\left(r_1,r_2\right):= B\left(0,r_2\right)\setminus B\left(0,r_1\right) .
\]

Sometimes (and it will be clear from the context), the area of this annulus will be also denoted
$A\left(r_1,r_2\right)$.
\begin{lem}
\label{linearbehavior}
Fix $k>1$. Then the area of the annulus $A\left(r,kr\right)$ 
in the cigar soliton satisfies
\[
A\left(r,kr\right)\sim C_k r .
\]
\end{lem}

\bfseries \textit{Proof. }\normalfont The soliton metric
is given by
\[
ds^2=\frac{dx^2+dy^2}{1+x^2+y^2} .
\]

Then, if $P=\left(x,y\right)$ we have (we denote by $d$ the distance
function in the manifold)
\[
\begin{array}{rcl}
d\left(0,P\right)&=&\int_0^{\alpha}\frac{1}{\sqrt{1+\rho^2}}\,d\rho \quad \mbox{where}\quad \alpha=\sqrt{x^2+y^2}\\
&\sim& \log \alpha .
\end{array}
\]

Therefore, if $r=\log \alpha$,
\[
\begin{array}{rcl}
A\left(\log \alpha,k\log \alpha\right)&\sim &2\pi \int_{\alpha}^{\alpha^k} \frac{r}{1+r^2}\,dr\\
&=& \pi \log \left(1+\alpha^k\right)-\pi\log \left(1+\alpha^2\right)\\
&\sim& \frac{k}{2}\pi \log \alpha .
\end{array}
\]

This shows the claim.

\hfill $\Box$

We proceed now with the argument as suggested by Perelman in \cite{Per1}. Given $r>0$, define the following function
\begin{equation}
\label{testfunction1}
\phi =
\left\{
\begin{array}{l}
1 \quad \mbox{if}\quad p\in A\left(2r, 3r\right)\\
\frac{\epsilon}{r} \quad \mbox{if} \quad p\in B\left(0,r\right)\\
\epsilon\exp\left(-d\left(0,p\right)\right) \quad \mbox{if}\quad p\notin B\left(0,4r\right)  ,
\end{array}
\right.
\end{equation}
where $\epsilon>0$ is a very small number, and 
\begin{equation}
\left|\nabla \phi \right|\leq \frac {1}{r} .
\end{equation}

Finally, define
\begin{equation}
f = -\log \phi + c ,
\end{equation}
where $c$ is a constant such that
\begin{equation}
\frac{1}{4\pi r^2}\int_M \phi \exp\left(-c\right)\,dV=1 .
\end{equation}

In the case of the soliton metric, for $\phi$ defined as in (\ref{testfunction1})
we have
\[
\begin{array}{rcl}
\frac{1}{4\pi r^2}\int_M \phi \exp\left(-c\right)\,dV&\sim&
\int_{A\left(r,4r\right)}\exp\left(-c\right)\,dV\\
&\sim& \frac{1}{4\pi r}\exp\left(-c\right)\quad \mbox{by Lemma \ref{linearbehavior}}
\end{array}
\]
which implies
\[
c \sim -\log\left(r\right) .
\]

Also, we have the following well known estimate on the decay of the curvature on the soliton metric
\begin{prop}
The curvature satisfies the following estimate
\[
R\left(P\right)\sim \exp\left(-d\left(0,P\right)\right) .
\]
\end{prop}

\bfseries\textit{Proof. }\normalfont 
By the transformation law of the curvature under conformal change of metric, we have
for a point $P=\left(x,y\right)\in \mathbf{R}^2$,
\[
R\left(P\right)=\frac{4}{1+x^2+y^2} .
\]

But, as we computed in Proposition \ref{linearbehavior}, 
\[
\sqrt{x^2+y^2}\sim \exp\left(d\left(0,P\right)\right)
\]
and the statement of the Proposition follows.

\hfill $\Box$

We let $\tau=r^2$ in (\ref{perelman0}). Then,
\[
\begin{array}{rcl}
\int_M \tau\left|\nabla f\right|^2\exp\left(-f\right)\cdot \frac{1}{4\pi \tau}\,dV &=&
\int_M 4 r^2\left|\nabla \exp\left(-\frac{f}{2}\right)\right|^2\exp\left(-c\right)\cdot \frac{1}{4\pi r^2}\,dV\\
&\leq& \int_{A\left(r,4r\right)}4r^2 \frac{1}{r^2}\cdot \frac{\exp\left(-c\right)}{4\pi r^2}\,dV ,
\end{array}
\]
which remains bounded as $r\rightarrow \infty$. 

Since $R\leq \exp\left(-r\right)$, the integral 
\[
\int_M R \cdot r^2 \frac{1}{4\pi r^2}\exp\left(-f\right)\,dV
\]
remains bounded.

It is clear also that the expression
\[
\int_M n \frac{1}{4\pi r^2}\exp\left(-f\right)\, dV
\]
remains bounded.

Finally,
\[
\begin{array}{rcl}
\int_M f\exp\left(-f\right)\frac{1}{4\pi r^2}\,dV &=&
\int_M \phi\log \phi \frac{1}{4\pi r^2}\exp\left(-c\right)\,dV\\
&&-\int_M c \exp\left(-f\right)\cdot \frac{1}{4\pi r^2}\,dV .
\end{array}
\]

The first integral in the righthand side is bounded. Indeed,
\[
\begin{array}{rcl}
\int_M \phi\log \phi \frac{1}{4\pi r^2}\exp\left(-c\right)\,dV &\leq&
C\int_{A\left(r,4r\right)}\frac{1}{4\pi r^2}\exp\left(-c\right)\,dV\\
&\leq& C\frac{1}{r^2}\cdot r \cdot r .
\end{array}
\]

The second integral goes to $-\infty$ as $r\rightarrow \infty$.

Replacing all these estimates into (\ref{perelman0}), this choice of $\phi$ shows that,
if $g\left(t\right)$ is a solution to the Ricci flow, and if there is a 
blow up sequence $\left(p_l,t_l\right)$ converging to the soliton metric, the there is
a sequence of radii $r_l$ such that 
\[
\mu_l=\mu\left(g\left(t_l\right),r_l^2\right)\rightarrow -\infty \quad
\mbox{as} \quad l\rightarrow \infty .
\]

On the other hand, let $\tau = t_l-t+r_l^2$, and $\hat{f}_l\left(\cdot,t\right)$
be the solution in the interval $\left[0,t_l\right]$ to
\[
\left\{
\begin{array}{l}
f_t=-\Delta f+\left|\nabla f\right|^2-R+\frac{n}{2\tau}\\
\frac{\partial f}{\partial \nu}=0\\
f\left(\cdot,t_l\right)=f_l ,
\end{array}
\right.
\]
where $f_l$ is chosen so that $\mu\left(g\left(t_l\right),f_l,r^2_l\right)\leq \mu_l+1$.

Then, by the monotonicity formula, for $t=0$ 
\[
\mu\left(g\left(0\right),t_l+r_l^2\right)\leq
\mathcal{W}\left(g\left(0\right),\hat{f}_l\left(\cdot,0\right),t_l+r_l^2\right)
\leq \mu_l+1
\]
if $\lim_{l\rightarrow\infty}t_l=T$, it is not difficult to show that (see \cite{Rot1},\cite{Rot2}),
$\mu\left(g(0),t_l\right)\rightarrow \mu\left(g(0),T\right)$, and this gives a contradiction. This
shows the following

\vspace{.1in}
\noindent
{\bf Theorem.}
If $R_0>0$, and $k_{g_0}\geq 0$, then the Ricci flow as considered in this paper uniformizes the curvature
along a sequence of times.

\subsection{Remarks on the Proof of Theorem \ref{convexboundarycase}.}

Some of the results in this section can be improved. 

\begin{prop}
If $R>0$, 
$\int_M R \,dV +\int_{\partial M} k \,d\sigma > 0 $ and (for the unnormalized flow)
$\frac{R_{max}\left(t_i\right)}{R_{min}\left(t_i\right)}\rightarrow 1$ along a sequence of times $t_i\rightarrow T<\infty$ then
the normalized flow exists for all time.
\end{prop}

{\bf \textit{Proof. }} Since $\frac{R_{max}\left(t_i\right)}{R_{min}\left(t_i\right)}\rightarrow 1$, the area 
 of the surface (in the unnormalized flow) is going to 0
as $t_i\rightarrow T$ (by Bonnet-Myers). 
Gauss Bonnet and the fact that
we keep the geodesic curvature bounded, imply that 
\[
A\left(t\right) \sim C\left(T-t\right)
\]
which in turn implies that the normalized flow exists for all time.

\hfill $\Box$

Taking into consideration the results of this section,
the previous proposition shows that Theorem \ref{longtimeexistence} is only interesting in the case that $k_0\leq 0$.
 
Also, we want to substantiate the claim in the last paragraph of the introduction, we first show that,

\begin{prop}
\label{curvatureisbounded}
If the normalized flow exists for all time, then there is a constant $C>0$ such that 
$R<C$ for all time.
\end{prop}

\textit{Proof. } If it is not the case, by passing to the unnormalized flow, we can find
a sequence of times $t_i\rightarrow T$, where $T$ is the blow up time, such that
\[
\left(T-t_i\right)R_{max}\left(t_i\right)\rightarrow \infty \quad
\mbox{as} \quad i\rightarrow \infty
\]

This implies that there is a blow up limit which is 
an eternal solution to the Ricci flow (with either no boundary or
a totally geodesic boundary). But then this limit must be the cigar, which
in turn contradicts estimate of Theorem \ref{monotonicity2}.

\hfill $\Box$

Consider the normalized flow again. Let
\[
r=\frac{\int_M R\,dV}{Vol\left(M\right)}.
\]

Since
\[
\frac{d r}{dt}=\int_{\partial M}\left[k_g\left(R-r\right)-2k_g'\right] \,d\sigma
\]
and it can be shown that $l_t(\partial M)<C$ and $k\sim C\exp(-\delta t)$, one can prove that 
there is $\epsilon>0$ such that $r\geq \epsilon$. This shows that there is a constant $c>0$ such that
\[
R_{max}\left(t\right)\geq c \quad \mbox{where}\quad  
R_{max}\left(t\right)=\max_{x\in M}R\left(x,t\right) .
\]

Passing to the unnormalized flow we conclude that
$R_{max}\left(t\right)$ is comparable to $\frac{1}{T-t}$ where $T$ is the time of blow up. 
Therefore, given any sequence $t_k\rightarrow T$, a blow up argument will show that along the 
chosen sequence
\[
\lim_{k\rightarrow \infty}\frac{R_{max}\left(t_k\right)}{R_{min}\left(t_k\right)}=1 .
\]

Passing back to the normalized flow we obtain,

\begin{thm}
Assume $R>0$ and $k_0>0$ at $t=0$. Then for the normalized Ricci flow the following holds
\[
R_{max}\left(t\right)-R_{min}\left(t\right)\rightarrow 0 \quad \mbox{as} \quad t\rightarrow \infty .
\]
\end{thm}

Exponential convergence in our case of study should be expected. If one for example tries to apply the methods in
Section 6 of
\cite{Struwe} to show exponential convergence, the arguments therein should work well since the boundary terms
are decaying to 0 exponentially fast.

\section{Subsequential Convergence in the case $k_0<0$: A Family of Examples.}\label{FamiliesOfExamples} 

Again, having Theorem \ref{limitingbehavior} in our toolbox, one way to proceed is as follows: we
study a scaling invariant geometric property of the initial condition that is preserved by the Ricci flow,
and show that this property does not hold for the cigar. This is what we do in what follows.

Consider a rotationally symmetric metric
\[
ds^2=dr^2+f\left(r\right)^2d\omega 
\]
on $\mathbf{B}^2$ that satisfies the following properties:
\newline
(P1) $R>0$ and $R$ is radially decreasing;
\newline
(P2) $\frac{\partial R}{\partial r}=k_g R$ in $\partial \mathbf{B}^2$.

For this type of metrics we are going to show that the Ricci flow uniformizes
curvature in the sense described above. 

Consider the quantity $\inf \frac{l^2\left(\partial B_{\rho}\right)}{A\left(B_{\rho}\right)}$ 
over balls centered at the pole in a rotationally symmetric
manifold, where by $A\left(B_{\rho}\right)$ we denote the area of such ball and by 
$l\left(\partial B_{\rho}\right)$ the length of its
boundary. Then we have,

\begin{lem}
\label{isoperimetricineq}
For $R\geq 0$ we have that
\[
\inf_{\rho>0} \frac{l^2\left(\partial B_{\rho}\right)}{A\left(B_{\rho}\right)}=\frac{l\left(\partial M\right)^2}{A\left(M\right)} .
\]
\end{lem}
\bfseries \textit{Proof. }\normalfont 
Since 
\[
\frac{l^2\left(\partial B_{\rho}\right)}{A\left(B_{\rho}\right)}=
\frac{4\pi^2 f^2\left(\rho\right)}{2\pi \int_0^\rho f\left(r\right)\,dr} ,
\]
all we must show is that for $\rho>0$ the expression
\[
\frac{f\left(\rho\right)f'\left(\rho\right)}{f\left(\rho\right)}=f'\left(\rho\right)
\]
is nonincreasing. But this is not difficult to see since $f''=-Kf$, where $K$ is the
Gaussian curvature of $M$.

\hfill $\Box$

Condition (P2) guarantees that the solution to (\ref{BVP2}) is at least $C^3$, and hence
we can use the Maximum Principle to show the following,
\begin{prop}
\label{nonincreasing}
For metrics satysfying properties (P1) and (P2),
under the Ricci flow the scalar curvature $R$ remains radially nonincreasing.
\end{prop} 

\bfseries \textit{Proof. }\normalfont 
The solution metric can be written as
\[
ds^2=h\left(r,t\right)^2dr^2+f(r,t)^2 d\omega^2 .
\]

Let $\frac{\partial}{\partial \rho}$ be the unit vector 
in the direction of $\frac{\partial}{\partial r}$, i.e., $\frac{\partial}{\partial \rho}=\frac{1}{h}\frac{\partial}{\partial r}$
. Define $w=fR_{\rho}$. Then, from the evolution equation for the curvature, we have
for $\rho>0$
\[
R_t=\frac{1}{f}w_{\rho}+R^2 .
\]

Differentiate with respect to $\rho$ to get,
\[
R_{t\rho}=-\frac{f_{\rho}}{f^2}w_{\rho}+w_{\rho\rho}+2RR_{\rho} .
\]

On the other hand 
\[
R_{\rho t}=-\frac{1}{2}RR_{\rho}+R_{t\rho} ,
\]
which in turn implies,
\[
f\left(\frac{1}{2}RR_{\rho}+R_{\rho t}\right)=-\frac{f_{\rho}}{f}w_{\rho}+\frac{1}{f}w_{\rho\rho}+2Rw .
\]

Also,
\[
\left(fR_{\rho}\right)_t=f'R_{\rho}+fR_{\rho t}
\]

From the Ricci flow equation follows that,
\[
f'=-\frac{1}{2}Rf .
\]

This yields,
\[
fR_{\rho t}= \left(fR_{\rho}\right)_t+\frac{1}{2}RfR_{\rho} ,
\]
from where we obtain,
\[
\frac{1}{2}Rw+\left(fR_{\rho}\right)_t+\frac{1}{2}Rw=-\frac{f_{\rho}}{f}w_{\rho}+\frac{1}{f}w_{\rho\rho}+2Rw .
\]

This last equation simplifies to
\[
w_t=\frac{1}{f}w_{\rho\rho}-\frac{f_{\rho}}{f}w_{\rho}+Rw ,
\]
and the boundary conditions are
\[
\left\{
\begin{array}{l}
w=0 \quad \mbox{at}\quad \rho=0\\
w<0 \quad \mbox{at} \quad \rho=\mbox{radius of the surface}\\
w|_{t=0}\leq 0 \quad .
\end{array}
\right.
\]

The result then follows from the Maximum Principle.

\hfill $\Box$

Consider the normalized flow now. From the previous proposition follows that 
\[
r \geq R_{min}\left(t\right)=R|_{\partial M \times \left\{t\right\}} ,
\]
and hence,
\[
l_t\left(\partial M\right)=\exp\left[\int_0^t \left(r-R\right)\,d\tau\right]l_0\left(\partial M\right)
\geq l_0\left(\partial M\right) .
\]

This last fact, together with Lemma \ref{isoperimetricineq}, shows that
\[
\inf \frac{l^2\left(\partial B_{\rho}\right)}{A\left(B_{\rho}\right)}=l_t\left(\partial M\right)\geq l_0\left(M\right) ,
\]
i.e., the isoperimetric ratio remains bounded away from 0 throughout the flow, and for the
cigar we have that this isoperimetric ratio is equal to 0. This shows that
we are in alternative (A) of Theorem \ref{limitingbehavior}, and hence the
flow is uniformizing the curvature along a subsequence of times.

\subsection{A family of examples. } 

Now we want to construct explicitly a family of metrics satisfying conditions (P1) and (P2). 
Consider the family of metrics $ds^2=dr^2+f_{\epsilon}\left(r\right)^2 d\omega^2$, where
\begin{equation}
\label{0}
f_{\epsilon}(r)=\left(1-\epsilon\right)\sin r + \epsilon r.
\end{equation}

Then, an elementary calculation yields
\begin{eqnarray*}
f'=\left(1-\epsilon\right)\cos r +\epsilon\\
f''=-\left(1-\epsilon\right)\sin r .
\end{eqnarray*}

The Gaussian curvature is given by
\[
-\frac{f''}{f}=\frac{\left(1-\epsilon\right)\sin r}{\left(1-\epsilon\right)\sin r +\epsilon r} ,
\]
which can be rewritten as
\[
K=\frac{1}{1+\left(\frac{\epsilon}{1-\epsilon}\right)\frac{r}{\sin r}} .
\]

This expression shows that $K$ is decreasing in $r$
(in fact, $\frac{r}{\sin r}$ is increasing). On the other hand, the
expression for the mean curvature of the geodesic spheres is given by
\[
\frac{f'}{f}=\frac{\left(1-\epsilon\right)\cos r+\epsilon}{\left(1-\epsilon\right)\sin r+\epsilon r} .
\]

The equation we want to be satisfied for certain value of $r$ is 
\[
f'''=2\frac{f''f'}{f} ,
\]
that is to say,
\[
\left(1-\epsilon\right)\cos r=\frac{2\left(1-\epsilon\right)\sin r\left[\left(1-\epsilon\right)\cos r+\epsilon\right]}
{\left(1-\epsilon\right)\sin r+\epsilon r} ,
\]
which simplifies to
\begin{equation}
\label{1}
\epsilon r \cos r-2\epsilon \sin r-\left(1-\epsilon\right)\cos r \sin r=0 .
\end{equation}

This is the equation to be solved. First of all notice that the expression in the left hand side
is always negative at $\frac{\pi}{2}$. At $r=\frac{3\pi}{4}$, the left hand side evaluates to
\[
\epsilon\left(\frac{3\sqrt{2}\pi}{8}+\frac{1}{2}-\sqrt{2}\right)+\frac{1}{2} ,
\] 
which is positive for any $\epsilon>0$. Therefore, there is 
$r_0\in\left(\frac{\pi}{2},\frac{3\pi}{4}\right)$
solving equation (\ref{1}). Notice that $f'(r_0)<0$: this
follows readily from the fact that $\frac{\partial R}{\partial r}<0$, $R>0$ and
$\frac{\partial R}{\partial r}=kR$.


\end{document}